\documentclass[10pt]{article}
\usepackage[tbtags]{amsmath}
\usepackage{epsfig,amstext,amssymb,amsthm,latexsym}
\allowdisplaybreaks[4]
\usepackage{amssymb,color}

\definecolor{c20}{rgb}{0.,0.7,0.}
\definecolor{c30}{rgb}{0.,0.,1.}
\definecolor{c40}{rgb}{1,0.1,0.7}
\definecolor{c50}{rgb}{1,0,0}

\def\cgrr#1{\textcolor{c20}{#1}}

\newtheorem{theo}{Theorem}[section]
\newtheorem{sat}[theo]{Proposition}
\newtheorem{de}[theo]{Definition}
\newtheorem{lem}[theo]{Lemma}
\newtheorem{exxa}[theo]{Example}
\newtheorem{korr}[theo]{Corollary}
\newtheorem{remark}[theo]{Remark}
\newtheorem{remarks}[theo]{Remarks}

\newcommand{\neprop}[1]{{Proposition \ref{#1}}}
\newcommand{\netheo}[1]{{Theorem \ref{#1}}}

\newcommand{\prooftheo}[1]{ \textsc{Proof of Theorem} \ref{#1} }
\newcommand{\proofprop}[1]{\textsc{Proof of Proposition} \ref{#1}}

\newcommand{\kb}[1]{\boldsymbol{#1}}
\newcommand{\vk}[1]{\kb{#1}}

\def\kal#1{{\cal{ #1}}}

\def\fracl#1#2{\biggr(\frac{#1}{#2} \biggl) }

\newcommand{\abs}[1]{\lvert #1 \rvert}

\newcommand{\norm}[1]{\lVert #1 \rVert}
\newcommand{\normS}[1]{\lVert #1 \rVert}

\newcommand{\pk}[1]{\mbox{\rm$\vk{P}$} \{#1\} }

\newcommand{\pb}[1]{\mbox{\rm$\vk{P}$}\Bigl \{#1 \Bigr \}}

\newcommand{\R}{\!I\!\!R}

\newcommand{\inr}{\in \R}

\newcommand{\ldot}{,\ldots,}

\newcommand{\limit}[1]{\lim_{#1 \to   \infty}}
\newcommand{\todis}{\stackrel{d}{\to}}
\newcommand{\ntoi}{n \to \infty }

\newcommand{\BQN}{\begin{eqnarray}}
\newcommand{\EQN}{\end{eqnarray}}
\newcommand{\BQNY}{\begin{eqnarray*}}
\newcommand{\EQNY}{\end{eqnarray*}}

\newcommand{\BS}{\begin{sat}}
\newcommand{\ES}{\end{sat}}
\newcommand{\BT}{\begin{theo}}
\newcommand{\ET}{\end{theo}}
\newcommand{\BK}{\begin{korr}}
\newcommand{\EK}{\end{korr}}
\newcommand{\EQD}{\stackrel{d}{=}}

\newcommand{\BD}{\begin{de}}
\newcommand{\ED}{\end{de}}

\newcommand{\BIT}{\begin{itemize}}
\newcommand{\EIT}{\end{itemize}}
\newcommand{\BDI}{\begin{description}}
\newcommand{\EDI}{\end{description}}

\newcommand{\BRM}{\begin{remarks}}
\newcommand{\ERM}{\end{remarks}}

\newcommand{\QED}{\hfill $\Box$}

\newcommand{\IF}{\infty}
\newcommand{\BTH}{\begin{theo}}
\newcommand{\ETH}{\end{theo}}
\newcommand{\BPR}{\begin{sat}}
\newcommand{\EPR}{\end{sat}}

\newcommand{\BEX}{\begin{exxa}}
\newcommand{\EEX}{\end{exxa}}

\newcommand{\BC}{\begin{cases}}
\newcommand{\EC}{\end{cases}}
\newcommand{\COM}[1]{}
\newcommand{\BL}{\begin{lem}}
\newcommand{\EL}{\end{lem}}

\newcommand{\thr}{\vk{t}}

\def\SI{\Sigma}
\def\SIJJ{\Sigma_{JJ}}
\def\SIII{\Sigma_{II}}
\def\SIIJ{\Sigma_{IJ}}
\def\SIJI{\Sigma_{JI}}
\def\SIM{\SI^{-1}}
\def\SIJJM{\Sigma_{JJ}^{-1}}
\def\SIIIM{\Sigma_{II}^{-1}}
\def\SIJIIM{\SIJI \SIIIM}
\def\thrI{\thr_{I}}
\def\thrJ{\thr_{J}}

\def\SIJJM{\SIJJ^{-1}}

\def\njd{ \{1 \ldot k\}}
\def\a{\vk{a}}
\def\b{\vk{b}}
\def\x{\vk{x}}
\def\y{\vk{y}}
\def\X{\vk{X}}

\def\U{\vk{U}}

\def\b{\thr}
\def\IB{I}
\def\JB{J}

\def\1d{\{1 \ldot d\}}
\def\njk{\njd}

\def\qp#1#2{\vk{Q}(#1,#2)}

\def\inrdoz{\R^k \setminus (-\infty, 0]^k}

\def\nouI{\normS{\thrI}}

\def\thr{\vk{t}_n}
\def\thrI{\vk{t}_{n,I}}
\def\thrJ{\vk{t}_{n,J}}
\def\tnn{\normS{\thrI}}
\def\wtn{ \beta_n}

\def\ntny{\alpha_n}

\def\prodtseiX{ \prod_{i\in I} \a_I ^\top \SIIIM \vk{e}_i}

\def\prodtasei{ \prod_{i\in I} \a_I ^\top \SIIIM \vk{e}_i}

\newcommand{\equaldis}{\stackrel{d}{=}}

\begin{document}

\begin{center}
\thispagestyle{empty}

{\Large Asymptotics for Kotz Type III Elliptical Distributions}\\

       \vskip 0.4 cm

         \centerline{\large Enkelejd Hashorva}

        \centerline{\textsl{Department of Mathematical  Statistics and Actuarial Science}}
        \centerline{\textsl{University of Bern, Sidlerstrasse 5}}
        \centerline{\textsl{CH-3012 Bern, Switzerland}}
       \centerline{\textsl{enkelejd.hashorva@Stat.Unibe.ch }}

\today{}

\end{center}

%%%%%%%%%%%%%%%%%%%%%%%%%%%%%%%%%%%%%%%%%%%%%%
{\bf Abstract:} Let $\X$ be a Kotz Type III elliptical random vector
in $\R^k, k\ge 2,$ and let $t_n,n\ge 1$ be positive constants such
that $\limit{n} t_n= \IF$. In this article we  obtain an asymptotic expansion of the tail probability $\pk{\X> t_n \vk{a}},
\vk{a}\inr^k$. As an application we derive an approximation for
the conditional excess distribution. Furthermore, we discuss the asymptotic dependence of
Kotz Type III triangular arrays and provide some details on the estimation of conditional
excess distributions and survivor function of Kotz Type III distributions.
\bigskip


\section{Introduction}
Consider a Kotz Type III elliptical random vector $\X$ in $\R^k,
k\ge 2,$ with the stochastic representation
\BQN \label{eq:A0}
\X\equaldis A^\top
R \vk{U},
\EQN
where $A\inr^{k\times k}$ is a non-singular matrix; the associated random radius $R>0$ has the tail asymptotic
behaviour
\BQN \label{eq:kotz:Fk} \pk{R> u}&=&  (1+o(1))p
u^{N}\exp(-q u^\delta ), \quad \text{   with } \delta>0,p>0,q>0, N\inr, \quad u\to \IF.
 \EQN
Furthermore, $R$ is independent of the random vector $\vk{U}$ which
is uniformly distributed over the unit $k-$sphere of $\R^k$. Here
$\equaldis$ denotes the equality of distribution functions and
$^\top$ is the transpose sign.

Prominent examples of the Kotz Type III elliptical random vectors
are the Gaussian ones where $R^2$ is chi-squared
distributed with $k$ degrees of freedom (see e.g., Kotz et al.\
(2000)) and the broader class of Kotz Type I elliptical random
vectors with  $R^2$ a Gamma distributed random variable. See Kotz (1975) and
Nadarajah (2003) for the main properties of Kotz Type I elliptical
random vectors.

Since $\X$ with stochastic representation \eqref{eq:A0} is an
elliptical random vector its basic distributional properties are
well-known; see e.g., Cambanis et al.\ (1981), Fang et al.\
(1990), or Kotz et al.\ (2000).

The main goal of this paper is to investigate some asymptotical
properties of the Kotz Type III elliptical random vectors. In the
recent papers Hashorva (2006a, 2007a,b) it is shown that such random
vectors have an asymptotic behaviour similar to that of the Gaussian
random vectors. This fact is at first sight surprising since we
do not specify the distribution of $R$, assuming only the asymptotic
relation in (1.2). The main reason for this similarity is the fact
that the associated random radius $R$ has distribution function in
the Gumbel max-domain of attraction.

Our primary interest in this paper is the tail asymptotic behaviour of
$\X$. Explicitly, if $t_n \vk{a},n\ge 1, \vk{a}\inr^k$ are given
thresholds in $\R^k, k\ge 2$ such that $\limit{n} t_n=\IF$ and
$\vk{a}$ has at least one positive component, then  $\limit{n}
\pk{\X> t_n \vk{a}} =0$. Of interest is therefore to determine the
speed of the convergence to 0 of the tail probability $\pk{\X> t_n
\vk{a}}$. The bivariate setup is discussed in Hashorva (2007a). We
extend that result to the multivariate setup by utilising some
general results for elliptical distributions obtained in Hashorva
(2007b). Two applications of the tail asymptotic expansion which we present here are: \\
a) an approximation of the conditional excess
distribution function, b) the asymptotic dependence of Kotz
Type III elliptical triangular arrays.\\
The Kotz Type distributions are encountered in various
statistical applications (see Nadarajah (2003)). In the light of our new asymptotic results it is possible to
estimate the survivor function and the conditional excess distribution of a Kotz Type III random vector.

Organisation of the paper: In the next section we present some preliminary results. The exact tail asymptotics
is discussed in Section 3 which is then followed by two short sections devoted to approximation of the conditional
excess distribution and the asymptotic independence of the Kotz Type III elliptical triangular arrays. In Section
6 we discuss briefly estimation of the survivor function and the conditional excess distribution.
Proofs of all the results are relegated to  Section 7.

\section{Preliminaries}
We shall introduce some standard notation. Let $\x=(x_1 \ldot x_k)^\top \inr^k$ be a vector in $\R^k,k\ge 2,$
and let in the following  $I$ be a non-empty index sets of $\njd$.
Denote by $\abs{I}$ the number of elements of $I$ and set $J:= \njd \setminus I$.  We define the subvector of $\x$ with respect to $I$
by $\x_I:=(x_i, i \in I)^\top \inr^{k}$. If $A\in \R^{k\times k}$ is a given matrix,
 then the submatrix  $A_{IJ}$ of $A$ is obtained by deleting both the rows and the columns of $A$
 with indices in $J$ and in $I$, respectively.  $A_{JI}, A_{JJ}, A_{II}$ are
similarly defined. We set $\SI:=A^\top A$ where $A$ is assumed to have a positive determinant $\abs{A}$ implying
that the inverse matrix $\SI^{-1}$ of $\SI$ exists. For notational simplicity we shall write $\x_I^\top, \SIJJM$ instead of $(\x_I)^\top, (\SI_{JJ})^{-1}$, respectively.
Given $\a, \x,\y\inr^k$ we shall define
\BQNY
%\x&\not=&\y, \text{ if } x_i \not= y_i, \text{ for at least one } i\le k, \\
 \x&>&\y, \text{ if } x_i>y_i,\quad  \forall\, i=1 \ldot k,\\
\x&\ge& \y, \text{ if } x_i\ge y_i, \quad \forall\,  i=1 \ldot k,\\
%\x &\not<  &\y, \text{ if } x_i \ge  y_i, \quad \text{ for some}\,  i=1 \ldot k,\\
\x+ \y &:=& (x_1+ y_1 \ldot x_k+ y_k)^\top,\\
c \x&:=&(cx_1\ldot cx_k)^\top,\quad c\inr,\\
\vk{a}\x&:=& (a_1 x_1 \ldot a_k x_k)^\top, \quad
\x/\vk{a}:= (x_1/a_1 \ldot x_k/a_k)^\top,\\
%$$a \x:=(ax_1\ldot ax_d)^\top,\quad a\inr, \quad
 \vk{0}&:=&(0\ldot 0)^\top \inr^k, \quad
 \vk{1}:=(1\ldot 1)^\top \inr^k, \text{  and  } \normS{\x_I}^2:= \x_I^\top \SIIIM\x_I.
\EQNY
Without loss of generality we shall assume that $\Sigma$ is a correlation matrix, i.e., all entries
of the main diagonal of $\SI$ are equal to 1.  If $\X$ is an elliptical random vector in $\R^k,k\ge 2$
with stochastic representation \eqref{eq:A0} in view of Lemma 12.1.2 in Berman (1992) we have
\BQN\label{eq:lema:Berman}
 X_i & \equaldis  & R U_1 , \quad 1 \le i  \le k,
 \EQN
 where $X_i$ is the $i$-th component of $\X$ and $U_1$ is the first component of $\U$.
If further the associated random radius $R$ has the tail asymptotics \eqref{eq:kotz:Fk}, then for any $x\inr$
\BQN
\label{eq:R:w}
\frac{ \pk{R> u+ x/ ( q \delta u^{\delta-1})}}{\pk{R> u}} %&=&
%(1+o(1))\Bigl( 1 + x/ ( q \delta u^{\delta})\Bigr)^N \exp\Bigl(-q u^\delta \Bigl[ ( 1 + x/ ( q \delta u ^{\delta}))^\delta- 1\Bigr]\Bigr)\notag\\
&\to& \exp(-x), \quad u\to \IF,
 \EQN
hence $R$ has distribution function $F$ in the max-domain of attraction of the
Gumbel distribution $\Lambda(x)=\exp(-\exp(-x)),x\inr$. From the extreme value theory
a distribution function $F$ with upper endpoint $\IF$ belongs to the max-domain of attraction of $\Lambda$ if
\BQN\label{eq:gumbel:m}
\limit{u} \frac{1- F(u+ x/w(u))}{1- F(u)}& = & \exp(-x), \quad \forall x\inr,
\EQN
where $w(\cdot)$ is a positive scaling function (see e.g, Resnick (1987),
Reiss (1989), Embrechts et al.\ (1997), Falk et al.\ (2004), Kotz and Nadarajah (2005) or de Haan and Ferreira (2006)).
If the associated random radius $R$ has tail asymptotics given by \eqref{eq:kotz:Fk},
then \eqref{eq:R:w} implies that the distribution function of $R$ is in the Gumbel max-domaain of attraction
with the scaling function $w(\cdot)$  defined by
\BQN \label{eq:ww}
 w(u) &=& (1+o(1))q \delta u^{\delta-1}, \quad u\to \infty.
 \EQN
When the associated random radius $R$ possesses the chi-squared distribution with $k$ degrees of freedom we obtain
\BQNY
\pk{R>u}%&=& \frac{1}{2^{k/2}\Gamma(k/2)}\int_{u^2}^\infty s^{k/2-1}\exp(-s/2) \,ds\\
&=&\frac{(1+o(1))\exp(-u^2/2) u^{k-2}}{2^{k/2-1}\Gamma(k/2)}, \quad
u\to \IF, \EQNY where  $\Gamma(\cdot)$ is the Gamma function. Hence
in this case \eqref{eq:kotz:Fk} holds with \BQN \label{eq:ne:Gauss}
p:=\frac{1}{2^{k/2-1}\Gamma(k/2)}, \quad q:=1/2, \quad \delta:=2,
\text{ and } N:=k-2 \EQN
implying that standard Gaussian random vectors belong to the class of the Kotz Type III elliptical random
vectors. Tail asymptotics of the Gaussian random vectors is discussed
Dai and Mukherjea (2001), Hashorva and H\"usler (2003), and Hashorva (2003, 2005) among several other papers.

The solution of the following quadratic programming problem:
\BQN \label{sol:x}
 \qp{\SI}{\vk{a}}: \text{minimise $ \normS{\x}^2$ under the linear constraint } \x \ge \vk{a},
 \EQN
with $\vk{a} \in \R^k\setminus (-\IF, 0]^k$ is the main ingredient for determining the tail asymptotics
under consideration.

\def\b{\vk{a}}
\def\tilb{\widetilde{\a}}

\BS \label{prop:pre:gaus}  %(Hashorva and Kotz (2007a))
Let $\Sigma \inr^{k \times k},k\ge 2$ be a positive definite
correlation matrix and let $\b\in \R^k \setminus (-\IF, 0]^k$ be a given vector. Then the quadratic programming problem
$ \qp{\SI}{\b}$ has a unique solution $\tilb$ defined by a   unique non-empty index set $\IB\subset \njd$ such that
\BQN \label{eq:IJ} \tilb_{\IB}=
\b_{\IB} >\vk{0}_I,    \quad \SIIIM \b_{\IB}&>&\vk{0}_I, \\
\label{eq:alfa} \min_{\x \ge  \b} \normS{\x}^2=\min_{\x \ge  \b}
\x^\top \SIM\x= \normS{\tilb}^2
= \normS{{\b_{\IB}}}^2 & = & \b_{\IB}^\top \SIIIM\b_{\IB}>0
\EQN
and in the case when $\abs{I}< k$ we have for $\JB :=\1d\setminus I$
\BQN
\label{eq:prop:A1} \tilb_{\JB}&=&   - ((\SI^{-1})_{JJ})^{-1}(\SIM)_{JI}\b_I=
\SIJI \SIIIM \b_{\IB}\ge \b_{\JB}.
\EQN
 Furthermore, for any $\x\inr^k$
\BQN \label{eq:new}
\x^\top \SIM \tilb= \x_I^\top
\SIIIM \tilb_I= \x_I^\top \SIIIM\b_I = \sum_{i\in I} x_i \vk{e}_i^\top \SIIIM \b_I,
\EQN
with $\vk{e}_i$ being the $i$-th unit vector in $\R^{\abs{I}}$ and $\vk{e}_i^\top \SIIIM \b_I>0, i\in I$. If $\b= c \vk{1}, c\in (0,\IF)$, we have $ 2 \le \abs{I} \le k$ where $\abs{I}$
denotes the number of elements of $I$.
\ES
\COM{
In view of the above proposition the minimum of  $\qp{\SI}{\vk{a}}$ is attained at $\vk{a}$ if
the Savage condition (Savage (1962))
 \BQN
 \label{eq:SAV}
\SIM \vk{a} &> &\vk{0}
\EQN
is  satisfied. Otherwise there exists a unique non-empty index set
$I \subset \{1 \ldot k \}$ which defines the unique solution of $\qp{\SI}{\vk{a}}$.
Below we will refer to the index set $I$
as the minimal index set.\\
Finally note that since
$$ \norm{c\x}= c\norm{\x}, \quad \forall c>0, \x\inr^k$$
this implies that the unique solution of $\qp{\SI}{t\vk{a}}, t>0,\vk{a}\inr^k$
coincides with the unique solution of $\qp{\SI}{\vk{a}}$ multiplied by $t$.
}
Below we will refer to the index set $I$
as the minimal index set. Note in passing that
$$ \norm{c\x}= c\norm{\x}, \quad \forall c>0, \quad \x\inr^k, $$
hence the unique solution of $\qp{\SI}{t\vk{a}}, t>0,\vk{a}\inr^k$
coincides with the unique solution of $\qp{\SI}{\vk{a}}$ multiplied by $t$.
 Next, we provide a general result for elliptical random vectors which follows from Theorem 3.1
and Theorem 3.4 in Hashorva (2007b). If $\vk{t}_n$ is a given vector in $\R^d$ then we write for notational simplicity
$\vk{t}_{n,K}$ instead $(\vk{t}_n)_K$ with $K \subset \njd$.

\def\tnnw{\tnn w(\tnn)}

\BT \label{theo:perturb} Let $\X\EQD  A^\top  R \U$ be an elliptical random vector in $\R^k, k\ge 2$ where $R>0$ has distribution function $F$ with an infinite upper endpoint being independent of $\vk{U}$ that is uniformly distributed on the unit sphere of $\R^k$ and $A$ is a non-singular $k$-dimensional square matrix such that $\SI:=A^\top A$ is a correlation matrix.  Assume that the distribution function $F$
is in the Gumbel max-domain of attraction with the positive scaling function $w$ and let $t_n,n\ge 1$
 be positive constants converging to $\IF$. Denote by $I$ the minimal index set of the quadratic programming problem $  \qp{\SI}{\vk{a}}, \vk{a}\in \inrdoz.$
 Assume, for simplicity, that $\norm{\vk{a}_I}=1$.  If $\thr,n\ge 1$ are given vectors in $\R^k$ such that
\BQN
\limit{n}  w(t_n) (\thrI- t_n\vk{a}_I)
&=& \vk{q}_I\inr^{m},
\EQN
where $m:=\abs{I}$ and furthermore if $m< k$ for $J:=\njd \setminus I$
\BQN
\limit{n}\sqrt{t_n w(t_n)}( \thrJ/t_n -    \SIJI\SIIIM \vk{a}_I)&=&
 \vk{q}_J\in [-\IF, \IF)^{k-m}
\EQN holds, then we have \BQN \label{eq:ttn} \pk{\X> \thr}&= &
(1+o(1))\exp(- \vk{q}_I^\top \SIIIM \vk{a}_I)
\frac{\Gamma(k/2)2^{k/2 -1} \pk{\vk{Z}_J > \vk{q}_J \lvert \vk{Z}_I=
\vk{0}_I}}{(2\pi)^{\abs{I}/2}
\abs{\SIII}^{1/2} \prodtseiX} \notag\\
&& \times (t_n w(t_n) ) ^{1-(k+m)/2} (1- F(t_n)), \quad \ntoi, \EQN
where $\vk{Z}$ is a standard Gaussian random vector in $\R^k$ with
covariance matrix $\SI$, $\pk{\vk{Z}_J
> \vk{q}_J \lvert \vk{Z}_I= \vk{0}_I}=1$ if $m=k$, and $\vk{e}_i$ is
the $i$-th unit vector in $\R^{m}$.
\ET

\def\b{\vk{a}}
\def\tilb{\widetilde{\a}}

\section{Exact Tail Asymptotics}
As shown in \netheo{theo:perturb} the exact asymptotic behaviour of elliptical random vectors can be derived provided that
the distribution function of the associated random radius $R$ is in the Gumbel max-domain of attraction.
\def\nouI{\ntny}
\def\tnn{\ntny}
\def\tnnw{\tnn \wtn}
\def\by{\norm{\a_I}}
As mentioned previously the associated random radius of Kotz Type
III elliptical random vectors is in the Gumbel max-domain of
attraction. This fact and the above theorem lead us to the following
result:
%exact tail asymptotiSince by defiWe state now the main
%result of this section.

\BT\label{theo:main1} Let $\X$ be a Kotz Type III elliptical random vector in $\R^k, k\ge 2$ where the associated random radius $R$
has the tail asymptotic given in \eqref{eq:kotz:Fk}, and $A\inr^{k\times k}$ be a non-singular matrix
with $\SI:=A^\top A$ a correlation matrix. Let $t_n,n\ge 1$ be a positive sequence
and let $\vk{a}\inr^k \setminus (-\IF, 0]^k$ be a given vector. Denote by $I$ the minimal index set related to the quadratic programming problem $\qp{\SI}{\vk{a}}$
and define $\vk{v}_n\inr^k, n\ge 1$ by
$$(\vk{v}_n)_I:=%w(t_n \norm{\vk{a}_I})=
q \delta (t_n \norm{\vk{a}_I})^{\delta -1}\vk{1}_I , \quad
\text{ and for   $\abs{I}< k$ set } \quad (\vk{v}_n)_J:= \sqrt{q \delta} (t_n \norm{\vk{a}_I})^{\delta/2 -1}\vk{1}_J.$$
If $I$ has $m$ elements and $\limit{n}t_n= \IF$, then for any $\x\inr^k$ we have
\BQN
\label{eq:theo1:1}
\lefteqn{\pk{\X > t_n \vk{a}+ \x/\vk{v}_n}}\notag \\
&=& (1+o(1)) p( q \delta)^{(1-(k+m)/2)} {\by}^{ k+N+\delta (1-(k+m)/2)}
\frac{\Gamma(k/2)2^{k/2 -1}}{(2\pi)^{m/2}
\abs{\SIII}^{1/2} \prodtasei}  \exp(- \vk{x}_I^\top \SIIIM \vk{a}_I/\by ) \notag \\
&&\times \pk{\vk{Z}_J > \IF (\vk{a}_J - \SIJIIM  \vk{a}_I) + \x_J
\lvert \vk{Z}_I= \vk{0}_I} t_n ^{N+\delta (1-(k+m)/2)} \exp(- q\by^\delta t_n ^\delta), \quad n\to \IF,
\EQN
where $\vk{Z}$ is a standard Gaussian random vector in $\R^k$ with covariance matrix $\SI$.\\
Set $\pk{\vk{Z}_J > \IF (\vk{a}_J - \SIJIIM  \vk{a}_I) + \x_J
\lvert \vk{Z}_I= \vk{0}_I}$ to 1 if $J$ is empty and  substitute above $\IF\cdot 0$ by $0$.
\ET

\begin{remark}
a) In the above theorem the index sets $I,J$ do not depend on the
choice of $\x\inr^k$.

b) If $\vk{S}$ is a Kotz Type III random vector with covariance
matrix $\kal{I}$ which is the identity matrix we have by
the amalgamation property (see e.g., Fang et al.\ (1990)) of the
spherical random vector that
$$ \vk{S}\equaldis (I_1 \abs{S_1}\ldot I_k \abs{S_k}),$$
where $I_1 \ldot I_k$ are independent random variables taking values
$-1,1$ with probability $1/2$. The above theorem can be easily
extended to the case of weighted Kotz Type III elliptical random
vector defined via the stochastic representation
$$ \X\equaldis A (I_1^* \abs{S_1}\ldot I_k^* \abs{S_k})$$
where $A$ is a $k$-dimensional real square matrix and $I^*_i,i\le k$
are independent random variables taking two values with  $-1,1$,
with $\pk{I_i^*=1}\in (0,1], i\le k.$
\end{remark}

We present next an illustrating example.

{\bf Example 1}. Let $\X\inr^k, k\ge 2$ be a Kotz Type III random vector with underlying covariance matrix $\Sigma$ specified by
\begin{eqnarray}\label{eq:rho:1}
 \Sigma=(1-\rho )\kal{I} + \rho\vk{1}\vk{1}^\top, \quad  \rho \in (-1/(k-1 ), 1),
\end{eqnarray}
where $\kal{I}\in \R^{k\times k}$ is the identity matrix. We are interested on the asymptotic
behaviour of $\pk{\X > t_n \vk{1}}$ where $t_n,n\ge 1$ are  positive constants tending to infinity as $n\to \IF$.
Since necessarily $\rho> -1/(k-1)$ we have

\BQNY \SI^{-1} \vk{1} &=& \Bigl((1-\rho )\kal{I}\vk{1} + \rho(\vk{1}\vk{1}^\top)\Bigr)^{-1} \vk{1}
%&=& \Bigl(\frac{1}{1-\rho}\kal{I}- \frac{\rho}{(1-\rho)(1+ (k-1)\rho)} \vk{1} \vk{1}^ \top\Bigr) \vk{1}\\
%&=& \frac{1}{1-\rho}\vk{1}- \frac{\rho k}{(1-\rho)(1+ (k-1)\rho)} \vk{1} \\
= \frac{1}{1+ (k-1)\rho} \vk{1}> \vk{0},
\EQNY
hence the quadratic programming problem $\qp{\SI}{\vk{a}}$ has the unique
solution $t_n \vk{1}$ with minimal index set $I=\njk$. Simple calculations yield
\BQNY
\norm{\vk{1}}^2&=&  \vk{1}^\top \SI^{-1} \vk{1} =  %(1+ (k-1)\rho) \vk{1}^\top \vk{1} =
\frac{k}{1+ (k-1)\rho}=: C_{\rho}^2>0, \EQNY and \BQNY \abs{\SI}&=&
(1-\rho)^{k-1}(1+ (k-1)\rho) >0. \EQNY Assume for simplicity that
the parameters $q,\delta$ defining the tail asymptotic of the
associated random radius of $\X$ satisfy $q\delta=1$. Applying
\netheo{theo:main1} we obtain for any $\x\inr$ \BQNY
%\label{eq:theo1:1}
\lefteqn{ \pb{ \X > t_n \vk{1}+ (C_{\rho} t_n)^{1- \delta}\x }}\notag \\
%%%
%%%
\COM{ &=& (1+o(1)) p C_{\rho}^{ (k+N+\delta (1-k))/2}
\frac{\Gamma(k/2)2^{k/2 -1}}{(2\pi)^{k/2}
(1-\rho)^{k-1}(1+ (k-1)\rho) ^{1/2} (1+ (k-1)\rho)^{-k} } \\
&& \times  \exp \Bigl(- \vk{x}^\top \vk{1}/(1+ (k-1)\rho)\Bigr)\notag \\
&&\times
t_n ^{N+\delta (1-k)} \exp\Bigr(- q\fracl{kt_n}{1+ (k-1)\rho}^\delta \Bigr)\\
}
&=& (1+o(1)) p C_{\rho}^{ (k+N+\delta (1-k))/2}%\fracl{k}{1+ (k-1)\rho}^{ (k+N+\delta (1-k))/2}
\frac{\Gamma(k/2)2^{k/2 -1}(1+ (k-1)\rho)^{k-1/2} }{(2\pi)^{k/2}
(1-\rho)^{(k-1)/2} } \\
&& \times  \exp \Bigl(- \vk{x}^\top \vk{1}/\sqrt{k(1+
(k-1)\rho)}\Bigr) t_n ^{N+\delta (1-k)} \exp\Bigr(- q (C_{\rho}t_n)
^\delta\Bigr), \quad n\to \IF. \EQNY
\COM{ If $k=2$ we retrieve the result presented in Hashorva (2007a)
\BQNY \pb{ \X > t_n \vk{1}+  \Bigl
[t_n\sqrt{2/(1+ \rho)}\Bigr]  ^{1- \delta}\x } &=& (1+o(1)) p
\fracl{2}{1+ \rho}^{ 1+N/2-\delta/2 } \frac{(1+ \rho) ^{3/2}
}{(2\pi)
\sqrt{1-\rho}} \\
&& \times  \exp \Bigl(- (x_1+x_2)/\sqrt{2(1+ \rho)}\Bigr)
t_n ^{N-\delta } \exp\Bigr(- q\fracl{2}{1+ \rho}^{\delta/2} t_n^{\delta} \Bigr), \quad n\to \IF.
\EQNY
}
In the Gaussian case with $p,q,\delta$ as in \eqref{eq:ne:Gauss} we have %for $k\ge 2$
\BQNY
\lefteqn{ \pb{ \X > t_n \vk{1}+ ( C_{\rho} t_n) ^{1- \delta}\x }}\notag \\
%&=& (1+o(1))  \fracl{k}{1+ (k-1)\rho}^{ k+k-2+2(1-k)}
%\frac{1}{(2\pi)^{k/2}
%(1-\rho)^{k-1}(1+ (k-1)\rho) ^{1/2} (1+ (k-1)\rho)^{-k} } \\
%&& \times  \exp \Bigl(- \vk{x}^\top \vk{1}/(1+ (k-1)\rho)\Bigr)
%t_n ^{2+2(1-k)} \exp\Bigr(- \fracl{kt_n}{1+ (k-1)\rho}^\delta /2\Bigr)\\
&=& (1+o(1)) \frac{(1+ (k-1)\rho)^{k-1/2} }{(2\pi)^{k/2}
(1-\rho)^{(k-1)/2} } \exp \Bigl(- \vk{x}^\top \vk{1}/\sqrt{k(1+
(k-1)\rho)}\Bigr)\notag t_n ^{-k} \exp\Bigr(- (C_{\rho} t_n) ^2
/2\Bigr), \quad n\to \IF. \EQNY

\section{Approximation of the Conditional Excess Distribution}
In this section we consider an application of our asymptotic expansion derived in the previous section.
Consider  $\X$ a Kotz Type III random vector in $\R^k, k\ge 2$. In
various statistical applications it is of some interest to estimate the distribution function of
the excess random vector
$\X- t \vk{a}, t>0,\vk{a}\inr^k$ conditioning on the event $\X> t \vk{a}$. If $I,J$ are partitions of
$\njd$, it is of some interest to investigate the asymptotic behaviour of the distribution function of $\X_J$
given the partial information $\X_I > t \vk{a}_I$, with $t$ tending to infinity.\\
Thomas and Reiss (2007) provide a detailed treatment of statistical applications related to conditional distributions.
In the context of extreme value estimation of conditional distribution relies on asymptotic results derived for
large thresholds; in our context it means that further investigation of the quantities of interest requires that
$t$ tends to infinity. Explicitly, we discuss next the asymptotic properties $(n\to \IF)$ of two random
sequences $\vk{V}_{n;\vk{a}}, \vk{V}_{n;I;\vk{a}}, n\ge 1,\vk{a}
\inr^k$ defined via the stochastic representation
$$ \vk{V}_{n;\vk{a}}\equaldis  \X- t_n \vk{a} \lvert \X > t_n \vk{a}, \quad \text{ and }
\vk{V}_{n;I; \vk{a}}\equaldis \X_J \lvert \X_I > t_n \vk{a}_I, \quad
t_n\inr .$$
The asymptotic behaviour of both random
sequences under consideration are closely related to the tail
asymptotics of the Kotz Type III elliptical random vectors. We give next our first results.

\BT \label{kor:1}
Under the assumptions and the notation of \netheo{theo:main1} we have the convergence in distribution
\BQN
\vk{v}_n \vk{V}_{n;\vk{a}} & \todis \vk{W}, \quad n\to \IF,
\EQN
where $\vk{W}_I$ has independent components with survivor  function $\exp(- s\vk{a}_I^\top \SIIIM \vk{e}_i/\norm{\vk{a}_I}), i\in I, s>0$.
Furthermore if $\abs{I}< k$, then $\vk{W}_I$ is independent of $\vk{W}_J$ which has partial survivor function
\BQN\label{eq:w1}
 \pk{\vk{W}_L> \x_L}&=& \frac{\pk{\vk{Z}_{J} > \IF (\vk{a}_J - \SIJIIM  \vk{a}_I)+ \x^*_J \lvert \vk{Z}_I=
\vk{0}_I}}{\pk{\vk{Z}_J > \IF (\vk{a}_J - \SIJIIM  \vk{a}_I) \lvert \vk{Z}_I= \vk{0}_I}},
 \quad \forall \x\in [0,\IF)^{k}, \quad L \subset J
 \EQN
 where  $\x^*\inr^{\abs{J}}$ with $\x^*_L:=\vk{x}_L$
 and if $\abs{J}> \abs{L}$ set  $\x^*_{J\setminus L}:=\vk{0}_{J\setminus L}$.
\ET
As demonstrated in the main result of the previous section the
asymptotic expansion of the tail probability of interest is
determined via the unique index set $I$ of the related quadratic
programming problem. It is therefore easy to find  norming
constants $h_n, \vk{b}_n,n\ge 1$ so that  $h_n(\X- \vk{b}_n)_J$
given $\X_I> t_n  \vk{a}_I$ converges in distribution if $I$ is the
unique index set determining the solution of the quadratic
programming problem $\qp{\SI}{\vk{a}}$.
 We arrive thus at the following result:
\BT \label{kor:2} Under the assumptions and the notation of \netheo{theo:main1}, if further the index set $J$ is non-empty,
 then  we have the convergence in distribution
\BQN\label{eq:korr:2}
h_n  (\vk{V}_{n;I; \vk{a}}- t_n  \SIJIIM  \vk{a}_I   ) &
 \todis \vk{Z}_J \lvert \vk{Z}_I= \vk{0}_I, \quad n\to \IF,
\EQN
with $h_n:= \sqrt{q \delta} (t_n \norm{\vk{a}_I})^{\delta/2 -1}, n\ge 1$.
\ET
In view of \neprop{prop:pre:gaus} if for a given non-empty index set $I \subset \njk$
the vector $\a_I\inr^{\abs{I}}$ is such that $$\SIIIM \vk{a}_I> \vk{0}_I,$$
then the vector $\a^*\inr^k$ with components $\a^*_I:=\a_I, \a_J^*:=\SIJIIM  \vk{a}_I $ is the solution of the quadratic programming problem
$\qp{\SI}{\vk{a}^*}$. Consequently the above corollary can be formulated for every vector $\a$ and a non-empty index set
$I$ such that the vector $\SIIIM \vk{a}_I$ has positive components.

Instead of conditioning on the event $\X_I> t_n \vk{a}_I$ which is
initially dealt with in Berman (1982, 1983), since $\X_I$ possess an
absolutely  continuous distribution function when $\abs{I}< k$, we
may consider conditioning on $\X_I= t_n \vk{a}_I$. This was first
suggested in Hashorva (2006a). We reformulate Theorem 3.2 therein
for our specific setup. Define therefore a sequence of random vectors $\vk{V}_{n;I; \vk{a}}^*,n\ge 1$ in the
same probability space such that $\vk{V}_{n;I; \vk{a}}^*\equaldis \X_J \lvert \X_I = t_n \vk{a}_I, n\ge 1$.

\BT\label{eq:theo:E} Let $\X, \SI, t_n,n\ge 1$ be as in \netheo{theo:main1}.
If $I,J$ is a partition of $\njk$ and $\vk{a}\inr^k$ is such that $\norm{\a_I}>0$,
then we have the convergence in distribution
\BQN \label{eq:theo1:1}
h_n\Bigl( \vk{V}_{n;I; \vk{a}}^*- t_n \SIJI\SIIIM \vk{a}_I \Bigr) \lvert \X_I= t_n \vk{a}_J
& \todis & \vk{Z}_J \lvert \vk{Z}_I= \vk{0}_I, \quad n\to \IF,
\EQN
with $h_n:= \sqrt{q \delta} (t_n \norm{\vk{a}_I})^{\delta/2 -1}, n\ge 1$
and $\vk{Z}$ a standard Gaussian random vector in $\R^k$ with covariance matrix $\SI$.
\ET
\begin{remark}
a) The random vector $\vk{Z}_J \lvert \vk{Z}_I= \vk{0}_I$ is a Gaussian random vector
in $\R^{\abs{J}}$ with mean zero and positive definite covariance matrix $\SIJJ- \SIJI\SIIIM\SIIJ.$\\
b) It is remarkable that in both \eqref{eq:korr:2} and \eqref{eq:theo1:1} the same limiting random vector appears.
In \netheo{eq:theo:E} the index sets $I,J$ are not related to the quadratic programming problem
$\qp{\SI}{\vk{a}} $ which is the case in \netheo{kor:2}.
\end{remark}

{\bf Example 2}. Let $\X, C_{\rho}$ be as in Example 1 (recall we
set $q\delta=1$). Since the index set $J$ is empty we obtain
applying \netheo{kor:1} \BQN
 (t_nC_{\rho})^{\delta-1} \vk{V}_{n;\vk{1}} & \todis \vk{W},
\quad n\to \IF,
\EQN
where $\vk{W}$ has independent unit Exponential components.\\
In the Gaussian case  $\delta=2$ hence we have the convergence in distribution
\BQN t_nC_{\rho}\vk{V}_{n;\vk{1}} & \todis \vk{W}, \quad n\to \IF. \EQN

\section{Maxima of Triangular Arrays of Kotz Type III Random Vectors}
As the Gaussian distribution, the Kotz Type III multivariate distribution possess some interesting
asymptotic properties with respect to the asymptotic dependence and asymptotic behaviour of sample extremes.
In order to present those properties, we deal next with a random sequence of Kotz Type III vectors.\\
Let therefore $\X, \X_n, n\ge 1$ be independent random vectors in $\R^k,k\ge 2$ with common distribution $G$
such that $\X$ has stochastic representation (1.1). In view of \eqref{eq:lema:Berman} the marginal distributions of
$G$ are identical and furthermore Theorem 12.3.1 in Berman (1992) implies that each marginal distribution is
in the Gumbel max-domain of attraction with the scaling function $w(u)= q \delta u^{\delta- 1}, u>0$.
For any two components $X_i,X_j, i\not=j, i,j\le k$ we have applying \netheo{theo:main1}
$$ \limit{u}\frac{\pk{ X_i> u, X_j>u}}{\pk{X_i>u}}=0.$$
The above asymptotics implies that the sample maxima has independent components (see e.g., Reiss (1989)).
Asymptotic independence means that the componentwise sample maxima
$\vk{M}_n,n\ge 1$ converges to a random vector with independent unit
Gumbel components. Explicitly, we have
\BQN\label{huesler} \Bigl( (M_{n1}- b_n)/a_n \ldot (M_{nk}- b_n)/a_n \Bigr) \todis
(\kal{M}_1\ldot \kal{M}_k), \quad n\to \IF,
\EQN where
$$a_n:=b_n^{1- \delta}/(q \delta), \quad b_n:= G^{-1}_1(1- 1/n), \quad n>1,$$
and $G_1^{-1}$ is the inverse of the marginal distribution function $G_1$ of $G$.
H\"usler and Reiss (1989) have shown  that a triangular array of Gaussian random vectors
can be constructed such that the limiting distribution function of the sample maxima is a random vector in $\R^k$
with dependent components and max-stable multivariate distribution which possesses unit Gumbel marginal distribution.

In view of Hashorva (2006b) the same asymptotic results hold in the more general case of
the Kotz Type III distribution. Explicitly, let us consider the Kotz Type III multivariate elliptical
triangular array with stochastic representation \BQN\label{eq:ell}
(X_{n1}^{(j)} \ldot X_{nk}^{(j)}) ^\top \equaldis  A_n^\top R \vk{U}, \quad 1 \le j\le n, \quad n\ge 1,
\EQN
where $R>0$ has tail asymptotic behaviour as in \eqref{eq:kotz:Fk},
independent of $\vk{U}$ which is uniformly distributed on the unit sphere, and $A_n,n\ge 1$ is
a sequence of $k$-dimensional non-singular square matrix. If $\Sigma_n:=A_n^\top A_n, n\ge 1$ has
all main diagonal entries equal 1, then the convergence in distribution in \eqref{huesler} holds provided that
\BQN\label{eq:reiss}
\limit{n} ( \vk{1}\vk{1}^\top- \Sigma_n)  \frac{b_n}{2 a_n}
& =&  C \in (0,\infty)^{k\times k},
\EQN
with $a_n,b_n,n\ge 1$ given by
$$ a_n:= (q^{-1} \ln n )^{1/\delta -1}/(q \delta ),
\quad b_n:= (q^{-1}\ln n )^{1/\delta}+a_n\Bigl[ N\ln (q^{-1} \ln
n )/\delta+ \ln p \Bigr], \quad n> 1.
$$
If $q \delta=1, \delta=2,$ then we have
$$ a_n:= (2\ln n )^{ -1/2},
\quad b_n:= (2\ln n )^{1/2}+a_n\Bigl[ N\ln (2\ln
n )/2+ \ln p \Bigr], \quad n> 1.
$$
Condition \eqref{eq:reiss} in this case agrees with the one imposed in H\"usler and Reiss (1989) for the Gaussian
setup. The bivariate distribution $\kal{G}$ of $(\kal{M}_1, \kal{M}_2)$ is given by
$$ G_{\gamma}(x,y)= \exp \left( - \Phi \left(\gamma + \frac{x-y}{2\gamma}
\right)\exp(-y)- \Phi \left(\gamma+  \frac{y-x}{2\gamma}\right)\exp(-x)    \right ), \quad x,y\inr, $$
where
$\Phi$ is the standard Gaussian distribution function on $\R$ and $\gamma^2:= \limit{n} ( 1- \sigma_{12,n})  \ln n\in (0,\IF)$.

\section{Estimation of Joint Survivor and Conditional Excess Distribution}
Let $\X, \X_1 \ldot \X_n$ be independent random vectors in $\R^k,k\ge 2$ with common distribution function $G$ such that
$\X$ is a Kotz Type III random vector with stochastic representation (1.1) and matrix $A$ such that
$\Sigma:=A^\top A$ is a positive definite correlation matrix. In view of Lemma 12.1.2 in Berman (1992) the marginal distributions $G_i,i\le k $ of $G$ are equal. Furthermore, by Theorem 12.3.1 in Berman (1992) we obtain
\BQN
1- G_1(t)&=& (1+ o(1)) \frac{1}{2}\frac{\Gamma(k/2)}{\Gamma(1/2)} (t w(t))^{(k-1)/2} 2^{(k-1)/2}
\pk{R>t} , \quad t\to \IF.
\EQN
In various statistical applications given the finite sample $\X_1 \ldot \X_n, n>1,$ estimation of
the joint survivor probability and the conditional excess $\psi_t, \psi^*_{t,\x}$
$$\psi_t:= \pk{\X > t\vk{1}}, \quad \psi^*_{t,\x}:=\pk{\X- t\vk{1} > \x\lvert \X>t \vk{1}}, \quad \x\inr^k$$
 for $t$ large enough is of certain interest.

Let in the following $\hat \psi_{n;t}, \hat \psi_{n;t,\x}^*$ denote two estimators of these quantities
which we specify below.
In view of our asymptotic results in order to estimate $\psi_t$ we need to estimate $p,N$ and $q,\delta, \SI^{-1}$,
whereas for estimating $\psi^*_t$ we need to estimate only $p,\delta$ and $\SI^{-1}$.
We note in passing that if $\X$ is a Gaussian random vector, then $N=k-2$ and $\delta=2$.
 We  assume for simplicity below that in our setup of Kotz Type III elliptical random vectors
 the constants $p,N$ are known.

Estimation of $\SI$ is dealt with in several recent papers, see for instance see e.g.,
Schmid and Schmidt (2006), Schmidt and Schmieder (2007), or Sarr  and Gupta (2008). Let $\hat \SI_n^{-1}$ denote an estimator of $\SI_n^{-1}$.
\cgrr{We note that estimation of the precision matrix $\SI^{-1}$ is important, since we implicitly determine (estimate) the unique index set $I$
related to the quadratic programming problem $\qp{\SI}{\vk{1}}$. Under a more restrictive assumption on $R$, for instance $R^\alpha$ is
Gamma distributed with positive parameters $a,b$, then for estimating the dispersion matrix $\SI^{-1}$ we can utilise the recent results of
Sarr  and Gupta (2008).}

Estimation of $q$ and $\delta$ is closely related to the estimation of
the scaling function $w(u)= q\delta u^{\delta -1}, u>0$. As in Hashorva (2007c) we can estimate $q,\delta$ borrowing
the idea of Abdous et al.\ (2008).\\
Next, write
$Y_{1:n} \le \cdots \le Y_{n:n}$ for the associated order statistics of
$X_{i,1}, i\le n$ and define the following Gardes-Girard  estimator of $\delta$ by
$$ \hat \delta_n:= \frac{1}{T_n}\frac{1}{k_n}\sum_{i=1}^n \Bigl( \log Y_{n-i+1:n}-
\log Y_{n- k_n+1:n}\Bigr), \quad j=1,2,$$
with $1 \le k_n \le n, T_n>0, n\ge 1$ given constants satisfying
$$ \limit{n} k_n=\IF, \quad \limit{n} \frac{k_n}{n}=0, \quad
\limit{n} \log(T_n/k_n)= 1, \quad \limit{n} \sqrt{k_n} b(\log(n/k_n))\to \lambda\inr,$$
where $b$ is some  regularly varying function with index $-1$ related to the asymptotics of $\ln (1- G_1(t)),t\to \IF$
(see Gardes and Girard (2006)). The scaling coefficient $q$ can be estimated by (see Abdous et al.\ (2008))
\BQN
\hat q_{n}&:=& \frac{1}{k_n} \sum_{i=1}^{k_n}
\frac{ \log (n/i) }{(Y_{n-i+1:n})^{\hat \delta_n}}, \quad j=1,2, n>1.
\EQN
The estimators of $\hat \psi_{n;t}, \hat \psi_{n;t,\x}^*$ can now be defined by plugging in $\hat q_n,
\hat \delta_n, \hat \SI_n^{-1}$ in (3.16) and \eqref{eq:w1}, respectively. \cgrr{Based on the known asymptotic properties of these estimators
it is possible to construct further confidence intervals for both the survivor and the conditional excess function.}

\section{Proofs}
\proofprop{prop:pre:gaus} The claim follows from Proposition 2.1 in
Hashorva and H\"usler (2003)  and Proposition 2.1 in Hashorva (2005). \QED

\def\byy{h_n}
\def\by{\norm{\a_I}}
\prooftheo{theo:main1}
In view of \neprop{prop:pre:gaus} we have
$\vk{a}_I^\top \SIIIM \vk{e}_i> 0$ holds for all $i\in I$. Assume for simplicity that the index set
$I$ has less than $k$ elements. Then we have further $ \vk{a}_J\le  \SIJI \SIIIM \vk{a}_I, \norm{\a_J}>0.$
\eqref{eq:ww} implies that the associated random radius
$R$ has distribution function $F$ in the Gumbel max-domain of attraction with the
scaling function $w(u) = q \delta u^{\delta-1}, u>0$. Set
$$\byy:= t_n\norm{\a_I}>0, \quad \vk{t}_n:= t_n \vk{a}+ \x/\vk{v}_n, n\ge 1,$$
where $(\vk{v}_n)_I=w(\byy)\vk{1}_I$ and $(\vk{v}_n)_J=\sqrt{w(\byy)/\byy}\vk{1}_J$.  We have
$$ \limit{n} w(\byy) (\vk{t}_n -t_n \vk{a})_I = \x_I$$
and
\BQN\label{eq:wtn:thr} \limit{n} \fracl{ q \delta
\byy^{\delta-1} }{ \byy }^{1/2} t_n  \Bigl(\vk{a}_J - \SIJIIM
\vk{a}_I \Bigr) &=& \IF (\vk{a}_J - \SIJIIM  \vk{a}_I)\in [-\IF,
0]^{\abs{J}}, \EQN where we interpret $\IF\cdot 0$ as $0$. The
assumptions of \netheo{theo:perturb} are thus fulfilled, hence we
may further write
\BQNY \label{eq:ttn:b} \pk{\X> \thr}&= &
(1+o(1))\exp(- \vk{x}_I^\top \SIIIM \vk{a}_I/\norm{\a_J}) \by^{\abs{I}}\\
&& \times \frac{\Gamma(k/2)2^{k/2 -1} \pk{\vk{Z}_J > \IF (\vk{a}_J - \SIJIIM
\vk{a}_I) + \x_J\lvert \vk{Z}_I= \vk{0}_I}}{(2\pi)^{\abs{I}/2}
\abs{\SIII}^{1/2} \prodtasei}
(q \delta \byy  ^{\delta} ) ^{1+\abs{J}/2-k} \pk{R> \byy}\\
&= & (1+o(1))\exp(- \vk{x}_I^\top \SIIIM \vk{a}_I/\norm{\a_J}) ( q \delta)^{(1+\abs{J}/2-k)}\by^{\abs{I}}\\
&& \times \frac{\Gamma(k/2)2^{k/2 -1} \pk{\vk{Z}_J > \IF (\vk{a}_J - \SIJIIM  \vk{a}_I) +\x_J\lvert \vk{Z}_I= \vk{0}_I}}{(2\pi)^{\abs{I}/2}
\abs{\SIII}^{1/2} \prodtasei}
\byy  ^{\delta (1+\abs{J}/2-k)} \pk{R> \byy}
\\
&= & (1+o(1))\exp(- \vk{x}_I^\top \SIIIM \vk{a}_I/\norm{\a_J}) p( q \delta)^{(1+\abs{J}/2-k)} {\by}^{ \abs{I}+N+\delta (1+\abs{J}/2-k)}\\
&&\times
\frac{\Gamma(k/2)2^{k/2 -1} \pk{\vk{Z}_J > \IF (\vk{a}_J - \SIJIIM  \vk{a}_I) + \x_J\lvert \vk{Z}_I= \vk{0}_I}}{(2\pi)^{\abs{I}/2}
\abs{\SIII}^{1/2} \prodtasei} \notag\\
&& \times   t_n ^{N+\delta (1+\abs{J}/2-k)} \exp(- q \byy  ^\delta), \quad \ntoi,
\EQNY
hence the proof follows. \QED

\prooftheo{kor:1}
Let $L$ be a non-empty index set of $\njk$ and set $M:= \njk\setminus L$. For any $\x\inr^k$ we may write
\BQNY
\pk{\X> t_n \vk{a}}\pk{(\X- t_n \vk{a})_L> (\x /\vk{v}_n)_L\Bigl \lvert \X> t_n \vk{a}}
& =&\pk{ \X_L> t_n \vk{a}_L+ (\x/ \vk{v}_n)_L,  \X_M> t_n \vk{a}_M}\\
& =&\pk{\X> t_n \vk{a}+ \x^*/\vk{v}_n},
\EQNY
where $\x^*$ is a  vector in $\R^k$ with $\x^*_L:=\x_L, \x_M:=\vk{0}_M$. Assume for simplicity that
$\abs{I}< k$ (thus $J$ is non-empty). Utilising \netheo{theo:main1} for any $\x$ such that $\x_L\in [0,\IF)^{ \abs{L}}$ we have
\BQNY
\lefteqn{\limit{n} \pb{(\X- t_n \vk{a})_L> (\x /\vk{v}_n)_L\Bigl \lvert \X> t_n \vk{a}}}\\
& =&\limit{n} \frac{\pk{\X> t_n \vk{a}+ \x^*/\vk{v}_n}}{\pk{\X> t_n \vk{a}}}\\
&=& \exp(- \a_I^\top \SIIIM (\x^*)_I) \frac{\pk{\vk{Z}_{J} > \IF (\vk{a}_J - \SIJIIM  \vk{a}_I)+ \x^*_J \lvert \vk{Z}_I=
\vk{0}_I}}{\pk{\vk{Z}_J > \IF (\vk{a}_J - \SIJIIM  \vk{a}_I) \lvert \vk{Z}_I= \vk{0}_I}},
\EQNY
as $n\to \IF$, hence the proof follows. \QED

\prooftheo{kor:2}
For any non-empty index set $L\subset J$ set
$$K:=L\cup I, \quad w(u):= q \delta u^{\delta- 1}, u>0, \quad h_n:= \fracl{w(t_n\norm{\vk{a}_I})}{t_n\norm{\vk{a}_I}}^{1/2}>0,\quad n\ge 1.$$
We may write
\BQNY
\lefteqn{\pk{\X_I > t_n \vk{a}_I}\pk{h_n
(\X_J- t_n \vk{a}_I)_L> \y_L \Bigl \lvert \X_I > t_n \vk{a}_I}}\\
& =&\pk{ \X_L - t_n (\SIJIIM \vk{a}_I)_L> h_n^{-1} \y_L , \X_I > t_n \vk{a}_I}\\
& =&\pk{ \X_K>   t_n \vk{a}^*_K + h_n^{-1} \y^*_K}, \quad \y\inr,
\EQNY
where $\a^*,\y^*$ are vectors in $\R^{\abs{L}+\abs{I}}$ with
$$\a^*_L: = (\SIJIIM \a_I)_L, \quad \a^*_{I}:=\a_I, \quad \text{ and } \y^*_L:=\y_L, \quad \y^*_I:=\vk{0}_I.$$
By \neprop{prop:pre:gaus} $\SIJIIM \a_I \ge \a_J$, implying $\a^*_L  \ge \a_L$. Further since
$\SIIIM \vk{a}_I$ has all components positive we have that $\a^*$ is the unique solution of
the quadratic programming problem $\qp{B^{-1}}{\vk{a}^*}$, where $B:= \SI_{K,K}$.
Applying \netheo{theo:main1} %and using \eqref{eq:uniX}
we obtain as $n\to \IF$
\BQNY
\pk{h_n(\X_J- t_n \vk{a}_I)_L> \y_L \Bigl \lvert \X_I > t_n \vk{a}_I}& =&
\pk{\vk{Z}_{L} > \IF \vk{0}_L+ \y^*_L \lvert \vk{Z}_I=\vk{0}_I}\\
& =& \pk{\vk{Z}_{L} > \y_L \lvert \vk{Z}_I=\vk{0}_I},
\EQNY
hence the proof follows. \QED

%{\bf Acknowledgement:}  I thank Professor Samuel Kotz for several corrections and suggestions and a Referee.

\bibliographystyle{plain}

\end{document}